\documentclass[journal]{IEEEtran}

\usepackage{amsmath}
\usepackage{amsfonts}
\usepackage{graphics}
\usepackage{pstricks}
\usepackage{array}
\usepackage{float}
\usepackage{epsf}
\usepackage{graphicx}
\usepackage{multirow}
\usepackage{cite}
\usepackage{acronym}
\usepackage{comment}
\usepackage{upgreek}
\usepackage{subcaption}
\usepackage{nicefrac} 
\usepackage[font=small]{caption}

\newcommand{\bfg}[1]{\boldsymbol{#1}}

\newcommand{\e}[1]{\bfg e_{\rm #1}}

\usepackage{dotlessi}

\begin{document}

\title{Instantaneous Frequency Estimation in Unbalanced Systems Using Affine Differential Geometry}

\author{Ali Alshawabkeh, Georgios Tzounas,~\IEEEmembership{Member,~IEEE}, 
{\'A}ngel Molina-Garc{\'i}a,~\IEEEmembership{Senior Member,~IEEE}, 
\\
and Federico~Milano,~\IEEEmembership{Fellow,~IEEE}%
  \thanks{A.~Alshawabkeh, G.~Tzounas and F.~Milano are with the School of Electrical and Electronic Engineering, University College Dublin, Dublin, D04V1W8, Ireland.  
  {\'A}. Molina-García is with the Department of Electrical Engineering, Universidad Politécnica de Cartagena, Cartagena 30202, Spain. \\
  Corresponding author's e-mail: georgios.tzounas@ucd.ie}
  \thanks{This work is supported by the Sustainable Energy Authority of Ireland~(SEAI) by funding A.~Alshawabkeh and F.~Milano under project FRESLIPS, Grant No.~RDD/00681; and by the Seneca Foundation -- Science and Technology Agency of the Region of Murcia under the Regional Program for Mobility, Collaboration, and Knowledge Exchange ``Jim{\'e}nez de la Espada'' by funding {\'A}. Molina under Grant No.~22213/EE/23.}%
  \vspace{-5mm}
  }

\maketitle

\begin{abstract}
  The paper discusses the relationships between electrical and affine differential geometry quantities, establishing a link between frequency and time derivatives of voltage, through the utilization of affine geometric invariants.  Based on this link, a new instantaneous frequency estimation formula is proposed, which is particularly suited for unbalanced and single-phase systems.  Several examples as well as measurements based on two real-world events  illustrate the findings of the paper. 
\end{abstract}


\IEEEpeerreviewmaketitle

\vspace{-2mm}

\section{Introduction}
\label{sec.introduction}

The problem of frequency estimation has been studied for many years and several solution approaches have been reported, e.g.,~see \cite{liu2014three, song2022fast, reza2016accurate, nie2019detection, pradhan2005freq}.  
These approaches rely on a variety of methods, including phase-locked loops~(PLLs), discrete Fourier transform, Kalman filters, least squares, adaptive notch filters, etc.
Particularly for grid synchronization and control applications, PLLs are a popular solution due their performance and simplicity.  Three-phase PLLs, for example, are widely utilized to provide real-time phase/frequency estimations in grid-connected power converters. 
A conventional PLL configuration in three-phase system applications is the synchronous reference frame~(SRF) PLL, which relies on transforming input voltages to the $\rm dq$ synchronous reference frame and on regulating the frame's angular position so that either the $\rm d$- or $\rm q$-axis component is zero.  The analogous of SRF-PLL for single-phase systems is the quadrature signal generation~(QSG)-based PLL.  Given a single-phase voltage, the latter defines a second dimension through a fictitious quadrature signal, required to enable the application of the Park transform (and thus the formulation of $\rm dq$-axis voltage components), e.g.,~see \cite{hadjidemetriou2016synchronization}.    

Other approaches are based on the inverse Park \cite{santos2008comparison} and the Hilbert transform \cite{hao2007measuring}, and on second-order generalized integrators \cite{sahoo2021phase, 7562505}.  Although the above approaches provide robust frequency estimations under balanced conditions, they often perform poorly and result in estimations with sinusoidal ripple errors for unbalanced systems \cite{reza2019three, karimi2004estimation, meral2012improved, kulkarni2015design, golestan2012design}.  Reducing the bandwidth helps mitigate this issue and refine accuracy, but also compromises dynamic performance \cite{escobar2014cascade}.  Efforts to improve the performance of PLLs under unbalanced conditions include, among other studies, \cite{liu2014three, meral2012improved}.

PLLs belong to the broad family of time-domain methods.  In this paper we also focus on this family but approach the frequency estimation problem from an unconventional perspective, based on differential geometry.  The starting idea is that voltage vectors can be perceived as velocities of points on space curves and, as such, be analyzed using differential geometrical invariants. In our recent work, we defined these curves in a Euclidean space and, by applying the Frenet–Serret formulas, we derived a correspondence between curvature and instantaneous electrical frequency \cite{milano2021applications, milano2022paradox, milano2021geometrical, milano2022frenet}.  Despite providing accurate estimations for balanced systems, the curvature obtained in these works is time-varying in stationary unbalanced conditions, a result that clearly does not align well with the notion of angular frequency of stationary ac signals.

In this paper, we aim at solving this issue through an alternative theory of differential geometry of curves, namely through \textit{affine differential geometry}.  This theory has found applications in various areas, such as control of mechanical systems \cite{lewis2018bountiful}, computer vision \cite{craizer2006parabolic}, and motion identification \cite{flash2007affine}.  To the best of our knowledge, no application to power system analysis or frequency estimation has been proposed so far.  The reason for the utilization of affine geometry in this work is that affine geometry is intrinsically well suited to estimate the curvature of conic functions, e.g., ellipses and parabolas \cite{nomizu1994affine}.  As unbalanced conditions can be viewed as an elliptical curve of a three-phase voltage \cite{milano2022frenet}, affine geometry appears as an ideal approach to estimate the frequency.

The specific contributions of the paper are as follows.
\begin{itemize}
\item A derivation of the expressions for the affine arc length and curvature in terms of the voltage of an ac system.
\item A formula of the instantaneous frequency of a three-phase voltage as a function of affine geometric invariants.
\item A demonstration of the effectiveness of the proposed formula as a frequency estimation technique for unbalanced three-phase systems, as well as for single-phase systems.
\end{itemize}

The last two points are fully supported through a variety of examples.  The examples show in particular that, for unbalanced systems, the proposed expression yields a more precise instantaneous frequency estimation compared to PLLs and the Frenet-frame based method from \cite{milano2021applications}.

The remainder of the paper is organized as follows.  Section~\ref{sec.background} recalls basic concepts from affine geometry. These concepts are essential for the derivation of the theoretical results of the paper presented in Section~\ref{sec.theory}. Section~\ref{sec.case} tests the proposed approach through analytical and numerical examples. Finally, Section~\ref{sec.Conclusion} draws relevant conclusions.

\section{Outlines of Affine Differential Geometry}
\label{sec.background}

Affine geometry can be defined as a Euclidean geometry without measuring distances or angles \cite{flash2007affine}. 
Let us consider a smooth parametric curve in the plane:
\begin{equation}
  \bfg{x}(t) = x_1(t) \, \e{1} + x_2(t) \, \e{2} \, ,
\end{equation}
where $x_1(t), x_2(t) : \mathbb{R} \mapsto \mathbb{R}$ are smooth and $\e{1}$, $\e{2}$ form an orthogonal basis of the plane.
Let us also assume that $\bfg x$ does not have inflection points, i.e., the magnitude of the operator
\begin{equation}
  \label{eq:outer}
  [ \dot {\bfg x}(t) , \ddot {\bfg x}(t) ] \neq 0, \quad  \forall t \, ,
\end{equation}
never vanishes. In \eqref{eq:outer}, $\dot{\bfg x} = d\bfg{x}/dt$ and $\ddot{\bfg x} = d^2\bfg{x}/dt^2$, and the bracket operator $[ \bfg a , \bfg b ]$ of two vectors $\bfg a, \bfg b \in \mathbb{R}^2$, is $[ \bfg a , \bfg b ] = a_1 b_2 - b_1 a_2$.  The \textit{affine arc length} indicated with $\sigma$, is:
\begin{equation}
  \label{eq:sigma}
  \sigma(t) = \int_{t_0}^{t} [\dot {\bfg x}(t) , \ddot {\bfg x}(t)]^{1/3} dt \, ,
\end{equation}
or, equivalently:
\begin{equation}
  \label{eq: dot sigma}
  \dot \sigma(t)
  = {d \sigma(t)}/{dt}
  = [\dot{\bfg x}(t) , \ddot{\bfg x}(t)]^{1/3} \, .
\end{equation} 
A curve $\bfg{x}$ is said to be parameterized with $\sigma$ if, for all $\sigma$:
\begin{equation}
  \label{eq: cond}
  [ \bfg{x}'(\sigma) , \bfg{x}''(\sigma) ] = 1  \, ,
\end{equation}
where $\bfg{x}' = d \bfg x / d\sigma$ is the \textit{affine tangent} and $\bfg{x}'' = d^2 \bfg x / d\sigma^2$ is the \textit{affine normal}.  Applying the chain rule, $\bfg{x}'$ becomes:
\begin{equation}
  \label{eq: xprime}
  \bfg{x}'(\sigma) = \frac{d \bfg{x}}{d \sigma} =
  \frac{d \bfg{x}}{dt} \frac{dt}{d \sigma} = 
  \frac{\dot {\bfg{x}}(t)}{[\dot{\bfg x}(t) ,\ddot{\bfg x}(t)]^{1/3} } \, ,
\end{equation}
and, differentiating (\ref{eq: cond}) with respect to $\sigma$, one
obtains
$[\bfg{x}'(\sigma) , \bfg{x}'''(\sigma)] = 0$.
This result implies that $\bfg{x}'$ and $\bfg{x}'''$ are linearly
independent, leading to the relationship:
\begin{equation}
  \bfg{x}'''(\sigma) = -\kappa_a(\sigma) \, \bfg{x}'(\sigma) \, ,
\end{equation}
where $\kappa_a$ is the \textit{affine curvature} of $\bfg{x}$ and is defined as:
\begin{equation}
  \label{eq: affine curvature}
  \kappa_a(\sigma) = [\bfg{x}''(\sigma) , \bfg{x}'''(\sigma)] \, .
\end{equation}
The affine curvature is represented by the area of the parallelogram
formed by the vectors 
$\bfg{x}''$ and $\bfg{x}'''$.
It is relevant to note that for nonsingular conic sections, $\kappa_a$ is constant \cite{calabi1996affine}.  For $\kappa_a = 0$ the curve is a parabola; for $\kappa_a > 0$ an ellipse; and for $\kappa_a < 0$ a hyperbola.  In the next section, we consider $\kappa_a > 0$.

\section{Voltage in the Affine Plane}
\label{sec.theory}

The magnetic flux $\bfg \varphi$ is assumed to be the \textit{position} of a point on a space curve in generalized coordinates and, from Faraday's law, the \textit{speed} of such a point is the voltage \cite{milano2021geometrical}:
\begin{equation}
  \label{eq:flux}
  \bfg \varphi(t) \equiv -\bfg{x}(t) \, \quad \Rightarrow \quad 
  \bfg v(t) = - \dot {\bfg \varphi}(t) \equiv \dot{\bfg x}(t) \, .
\end{equation}

In \cite{milano2021applications}, it is shown that one can express voltage and current in terms of Frenet-frame coordinates and geometric invariants.  In the same vein, but using the coordinates and invariants of affine differential geometry, this section derives a new instantaneous frequency formula of electrical quantities.  We discuss only voltages, but the same procedure can be followed using currents.  We consider two scenarios, namely unbalanced three-phase systems; and single-phase systems.

\subsection{Three-Phase Unbalanced Voltages}

Let's assume that the phases $\rm abc$ of a three-phase voltage $\bfg v(t)$ constitute a set of orthogonal coordinates:
\begin{equation}
  \label{eq:vabc}
  \bfg v(t) = v_{\rm a}(t) \, \e{a} + v_{\rm b}(t) \, \e{b} + v_{\rm c}(t) \, \e{c} \, .
\end{equation}

The theory described in Section~\ref{sec.background} applies to curves in two dimensions. Thus, we first transform $\bfg v(t)$ into the shape:
\begin{equation}
  \label{eq:v12}
  \bfg{v}(t) = v_{1}(t) \, \e{1} +  v_{2}(t) \, \e{2} \, .
\end{equation}
This is conveniently achieved by applying the Clarke transform to \eqref{eq:vabc} and taking the $\alpha$ and $\beta$ components, as follows:
\begin{equation}
\label{eq:clarke}
  \begin{bmatrix}
    v_{\alpha}(t) \\
    v_{\beta}(t)
  \end{bmatrix}
  =
  \sqrt{\frac{2}{3}}
  \begin{bmatrix}
    1 & -\frac{1}{2} & -\frac{1}{2} \\
    0 & \frac{\sqrt{3}}{2} & -\frac{\sqrt{3}}{2} \\
  \end{bmatrix}
  \begin{bmatrix}
    v_{\rm a}(t) \\
    v_{\rm b}(t) \\
    v_{\rm c}(t)
  \end{bmatrix} .
\end{equation}
Thus, in \eqref{eq:v12} we have $v_1(t) = v_{\alpha}(t)$ and $v_2(t) = v_{\beta}(t)$.

\subsubsection{Stationary Sinusoidal Voltages}
\label{sub:ellipse}

We first consider an unbalanced stationary sinusoidal voltage, of which affine differential geometry allows obtaining the exact frequency.  Using Clarke's transform, the components of $\bfg{v}(t)$ in \eqref{eq:v12} are:
\begin{equation}
  \label{eq:basecase}
  v_1(t) = V_1 \cos \theta(t) \, , \qquad
  v_2(t) = V_2 \sin \theta(t) \, ,
\end{equation}
where $V_1$, $V_2$ are constant and $\theta(t) = \omega_o t + \theta_o$; $\omega_o$ is the fundamental synchronous reference frequency; $\theta_o$ is constant and its value depends on the chosen phase angle reference.

With the equivalence given in \eqref{eq:flux}, equation the time derivative of the affine arc length $\dot \sigma$ in \eqref{eq: dot sigma} can be written as:
\begin{equation}
  \label{eq:sigma_base}
  \dot \sigma
  = [ \bfg {v}(t) , \dot {\bfg v}(t) ] ^{1/3}
  = (\omega_o V_1 V_2)^{1/3} \, .
\end{equation}
Note that while $\bfg v$, $\dot {\bfg v}$ depend on time, $\dot \sigma$ does not.  Then, imposing that the voltage components are as in \eqref{eq:basecase}, one gets:
\begin{equation}
  \label{eq:xprime}
  \bfg{x}'(t)
  = {\bfg{v}(t)}/{\dot \sigma} \, ,  \
  \bfg{x}''(t)
  = {\dot {\bfg v}(t)}/{\dot \sigma^2} \, , \ 
  \bfg x'''(t) = {\ddot {\bfg v}(t)}/{\dot \sigma^3} \, ,
\end{equation}
where
\begin{equation}
  \label{eq:vdots}
  \begin{aligned}
    \bfg v(t) &= V_1 \cos\theta(t) \, \e{1} + V_2 \sin\theta(t) \, \e{2} \, , \\
    \dot{\bfg v}(t) &= - \omega_o V_1 \sin\theta(t) \, \e{1} +
    \omega_o V_2 \cos\theta(t) \, \e{2} \, , \\
    \ddot{\bfg v}(t) &= - \omega^2_o V_1 \cos\theta(t) \, \e{1} -
    \omega^2_o V_2 \sin\theta(t) \, \e{2} \, . \\
  \end{aligned}
\end{equation}
Then, using (\ref{eq: affine curvature}), (\ref{eq:sigma_base}), (\ref{eq:xprime}), the affine curvature $\kappa_a$ becomes:
\begin{equation}
  \label{eq:kappa_base}
  \kappa_a
  =  [\dot{\bfg v}(t) ,\ddot{\bfg v}(t)]/\dot \sigma^5
  = {\omega_o^3 V_1 V_2}/{\dot \sigma^5} \, ,
\end{equation}
where $\kappa_a$ is constant, which is as expected since \eqref{eq:basecase} describes an ellipse in the plane $(v_1, v_2)$.  Merging \eqref{eq:sigma_base}, \eqref{eq:kappa_base} we obtain:
\begin{equation}
  \label{eq:omega_o}
  \omega_o
  = \sqrt{\kappa_a} \, \dot{\sigma}
  = \sqrt{{[\dot{\bfg v}(t) ,\ddot{\bfg v}(t)]}/{[\bfg v(t) , \dot{\bfg v}(t)]}}\, ,
\end{equation}
which shows that, calculating the frequency of $\bfg v$, it suffices to measure it and estimate its first and second time derivatives.

\subsubsection{Transient Voltages}
\label{sub:vt}

In practice, noise, harmonics, and transients prevent deriving an explicit expression of the frequency.  Yet, in certain conditions, it is still possible to utilize the results above for a voltage of time-varying frequency and/or magnitudes.  Consider a time-varying voltage:
\begin{equation}
  \label{eq:vt}
  \bfg v (t) =
  V_1(t) \cos \vartheta (t) \, \e{1} +  
  V_2(t) \sin \vartheta (t) \, \e{2} \, ,
\end{equation}
where $\vartheta(t) = \omega_o t + \phi(t)$.  The conditions so that \eqref{eq:omega_o} holds for a voltage $\bfg v (t)$ in the form of \eqref{eq:vt} are:
\begin{align}
  \label{eq:cond1}
  \frac{d^h}{dt^h} \phi(t) &\ll \omega^h_o \, , \quad h = 1, 2 \, , \\
  \label{eq:cond2}
  \frac{d^h}{dt^h} \frac{V_i(t)}{\langle V_i \rangle}
                                &\ll \omega^h_o  \, , \quad i,h = 1, 2 \, ,
\end{align}
where $\langle \cdot \rangle$ denotes the average value.  For $h= 1$, \eqref{eq:cond1} indicates that the instantaneous frequency variation is close to the synchronous reference frequency; for $h=2$, \eqref{eq:cond1} imposes a boundary to the rate of change of frequency;

and \eqref{eq:cond2} impose that \textit{radial frequency} variations (see \cite{milano2021geometrical}) are small compared to the grid's fundamental frequency.  These assumptions are generally well satisfied in power systems.
Conditions \eqref{eq:cond1} and \eqref{eq:cond2} are sufficient for \eqref{eq:omega_o} to hold at least as a first order approximation.  In fact, the time derivative of $\bfg v (t)$ in \eqref{eq:vt} is:
\begin{equation*}
    \dot{\bfg v} = (\dot{V}_1 \cos \vartheta - \dot{\vartheta} V_1 \sin \vartheta) \, \e{1} 
    + (\dot{V}_2 \sin \vartheta + \dot{\vartheta} V_2 \cos \vartheta) \, \e{2} \, ,
\end{equation*}
and the second time derivative is:
\begin{equation*}
  \begin{aligned}
    \ddot{\bfg v} =
    &- \big (
    \ddot{\vartheta} V_1 \sin \vartheta +
    \dot{\vartheta}^2 V_1 \cos \vartheta  
    + \dot \vartheta \dot V_1 \sin \vartheta -
    \ddot{V}_1 \cos \vartheta
    \big ) \, \e{1} 
    \\&
    + \big (
    \ddot{\vartheta} V_2 \cos \vartheta +
    \dot{\vartheta}^2 V_2 \sin \vartheta 
    + \dot \vartheta \dot V_2 \cos \vartheta -
    \ddot{V}_2 \sin \vartheta
    \big ) \, \e{2} \, ,
  \end{aligned}
\end{equation*}
where the time dependency is omitted for economy of notation.

It is easy to show that by applying \eqref{eq:cond1} and \eqref{eq:cond2}, $\dot{\bfg v}$, $\ddot{\bfg v}$ can be approximated with the second and third equations of \eqref{eq:vdots} and, hence, the frequency can be estimated using \eqref{eq:omega_o}.
In summary, \eqref{eq:cond1}, \eqref{eq:cond2} lead to the following approximated expression of the instantaneous frequency of a time-varying unbalanced voltage:
\begin{equation}
  \label{eq:IF}
  \dot{\vartheta}(t)
    \approx \boxed{\omega_{a}(t) =
    \sqrt{ \frac{[\dot{\bfg v}(t) , \ddot{\bfg v}(t)]}{[\bfg v(t) , \dot{\bfg v}(t)]} } }
\end{equation}
The expression of $\omega_a$ in \eqref{eq:IF} is the main result of this work.  A fundamental condition for $\omega_a$ to work properly is that \eqref{eq:outer} is satisfied at all times.  Harmonics introduce points for which $[\dot{\bfg x}, \ddot{\bfg x}] \le 0$, which may lead to numerical issues.  However, these issues can be easily overcome if harmonics are adequately filtered.

\subsection{Application to Single-Phase Voltages}
\label{sub:single}

We consider a single-phase voltage with instantaneous value $v(t)$.  To apply the theory of Section~\ref{sec.background}, we first need to transform $v(t)$ into the shape of \eqref{eq:v12}. To this aim, we construct the second dimension by employing the voltage derivative, i.e.:
\begin{equation}
  \label{eq:v1ph1}
  v_1(t) = v(t) \, , \qquad v_2(t) = \dot{v}(t) \, .
\end{equation}
Since the time derivative of sinusoidal signals gives a $90^{\circ}$ rotation, using \eqref{eq:v1ph1} is equivalent to defining a quadrature axis.

\subsubsection{Stationary Sinusoidal Voltages}
\label{sub:1phase}

The result obtained in the previous section can be extended to a stationary sinusoidal single-phase voltage using \eqref{eq:v1ph1}.  Let the voltage be:
\begin{equation}
  v(t) = V \cos\theta(t) \, ,
\end{equation}
where $V$ is constant and $\theta$ is defined in \eqref{eq:basecase}.  Then, from \eqref{eq:v1ph1}, the components of the voltage vector are:
\begin{equation}
  \label{eq:basecase2}
  v_1(t) = V \cos\theta(t) \, , \qquad
  v_2(t) = -\omega_o V \sin\theta(t) \, .
\end{equation}
Substituting $V_1 = V$, $V_2 = \omega_o V$ in \eqref{eq:sigma_base}
and \eqref{eq:kappa_base}, one obtains:
\begin{equation}
  \label{eq:sigma_1ph}
  \dot \sigma
  = (\omega_o V)^{2/3} \, , \quad  \kappa_a
  = {\omega_o^4 V^2}/{\dot \sigma^5} \, .
\end{equation}
Apart from the fact that calculation of $\ddot{\bfg v}(t)$ in this case requires computing the third derivative of $v(t)$, equation \eqref{eq:omega_o} holds and allows estimating the frequency also for a single-phase voltage.

\subsubsection{Transient Voltages}
\label{sub:vt1ph}

Consider a time-varying voltage:
\begin{equation}
  \label{eq:v1t}
  v (t) = V(t) \cos\vartheta (t) \, ,
\end{equation}
where $\vartheta(t) = \omega_o t + \phi(t)$.  The voltage vector is defined as:
\begin{equation}
  \begin{aligned}
    \bfg v (t) = \;
    &V(t) \cos\vartheta (t) \, \e{1} \\ + \;
    & [\dot V(t) \cos\vartheta (t) -
      V(t) \dot \vartheta(t) \sin\vartheta(t) ] \, \e{2} \, .
  \end{aligned}
\end{equation}
If one assumes:
\begin{align}
  \label{eq:cond112}
  \frac{d^h}{dt^h} \phi(t) &\ll \omega^h_o \, , \quad 
  \ \
  \frac{d^h}{dt^h} \frac{V(t)}{\langle V \rangle} \ll \omega^h_o  \, , \quad \ \ 
  h = 1, 2, 3 \, ,
\end{align}
then \eqref{eq:IF} is also a good approximation of the instantaneous frequency of the time-varying single-phase voltage in \eqref{eq:v1t}.

\section{Case Studies}
\label{sec.case}

This section illustrates \eqref{eq:IF} in various conditions, comparing its accuracy with a SRF-PLL, as well as with the Frenet frame-based estimation $\omega_{\kappa}(t) = {[ \bfg v(t) , \dot{\bfg v}(t)]}/{|\bfg v(t)|^2}$, see \cite{milano2021applications}.  Voltage trajectories and frequency estimations are given in per unit (pu).  Formula \eqref{eq:IF} is calculated using a sampling of the voltage, applying the Clarke transform and then evaluating numerically the time derivatives of the $\alpha$, $\beta$ components.

\subsection{Three-Phase Voltage}

Let us consider the three-phase voltage vector given in \eqref{eq:vabc}:
\begin{equation}
  \bfg v(t) = v_{\rm a}(t) \, \e{a} + v_{\rm b}(t) \, \e{b} + v_{\rm c}(t) \, \e{c} \, ,
\end{equation}
with components:
\begin{equation}
  \begin{aligned}
    v_{\rm a}(t) &= V_{\rm a} \sin(\omega_o t + \phi_{\rm a}(t)) \, , \\ 
    v_{\rm b}(t) &= V_{\rm b} \sin(\omega_o t + \phi_{\rm b}(t) - \zeta_{\rm b}) \, ,  \\
    v_{\rm c}(t) &= V_{\rm c} \sin(\omega_o t + \phi_{\rm c}(t) + \zeta_{\rm c}) \, .  
  \end{aligned}
\end{equation}
Recall that we use \eqref{eq:clarke} to convert \eqref{eq:vabc} to the $(\alpha, \beta)$ plane.

\subsubsection{Balanced Voltage}

\begin{figure}[ht!]
  \centering
  \begin{subfigure}{0.49\columnwidth}
    \centering
    \includegraphics[width=\linewidth]{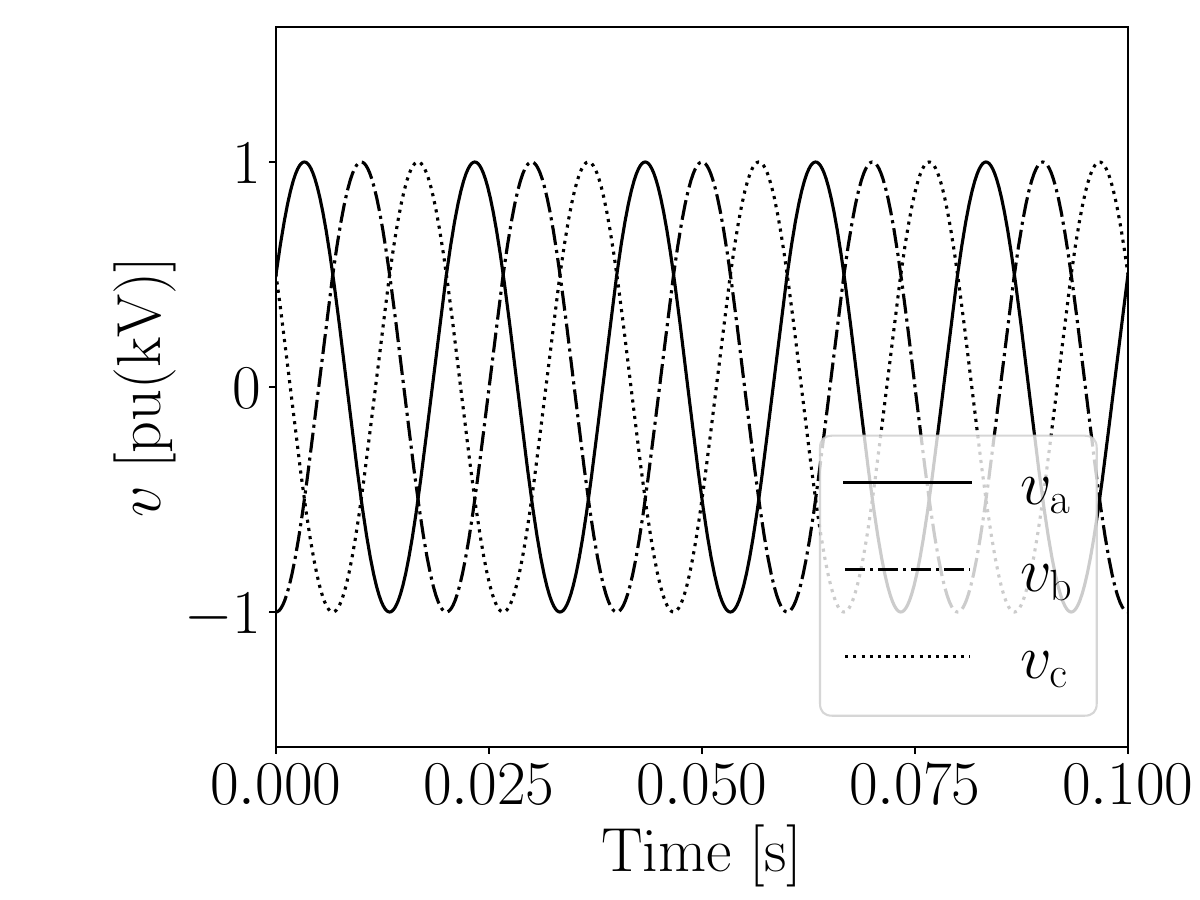}
    \vspace{-6mm}
    \caption{E1}
    \label{fig:E1a}
  \end{subfigure}
  \begin{subfigure}{0.49\columnwidth}
    \centering
    \includegraphics[width=\linewidth]{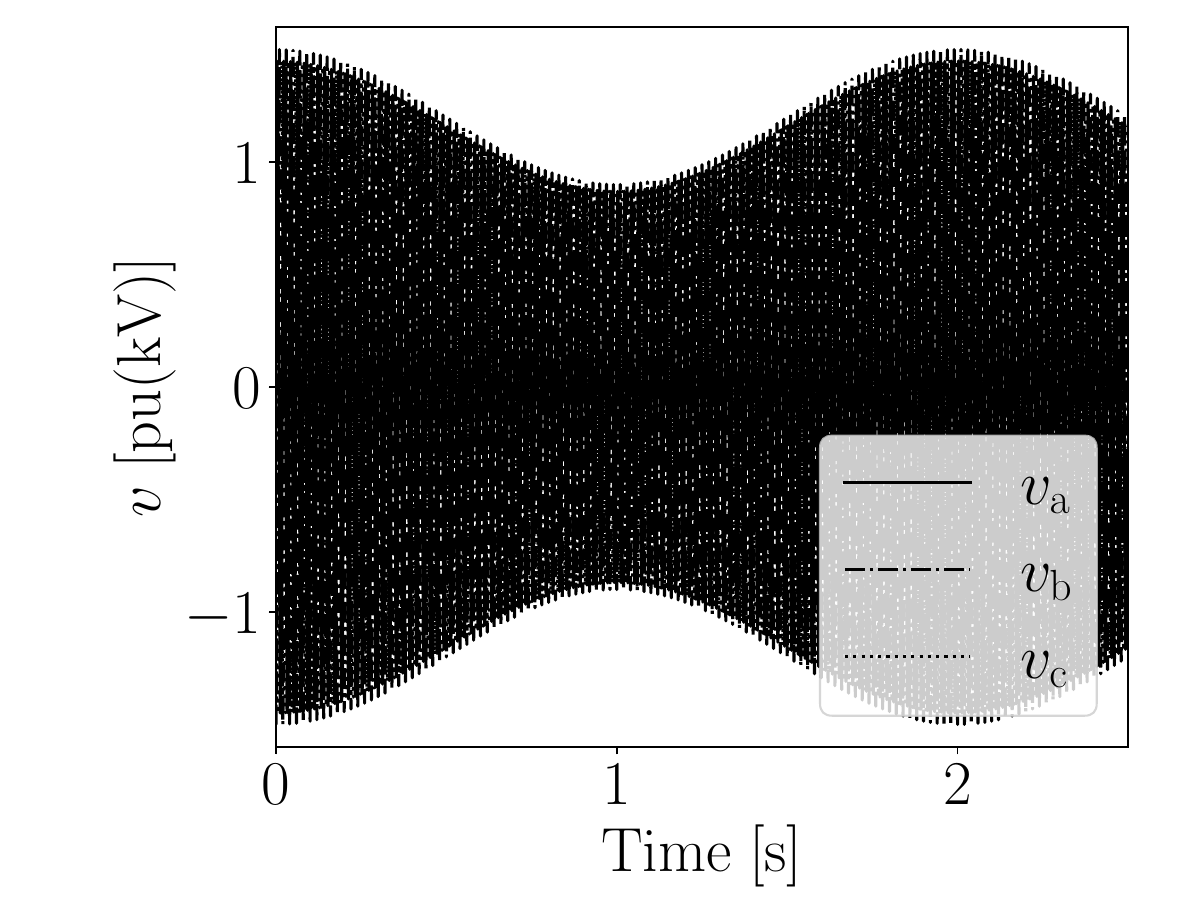}
    \vspace{-6mm}
    \caption{E2}
    \label{case:E2a}
  \end{subfigure}
  \vspace{-5.5mm}
  \caption{Balanced 3-phase voltage components.} 
  \label{fig:balanced system}
  \vspace{-5mm}
\end{figure}

We discuss two cases: (i) a stationary voltage, (ii) a voltage with time-varying magnitude.  In both cases, $\omega_o = 100 \pi $~rad/s.  The parameters used are:
\begin{itemize}
\item E1: $V_i= 12$~kV, $\phi_i = 0$ and $\zeta_{\rm b} = \zeta_{\rm c} = \frac{2 \pi}{3}$ rad.
\item E2: $V_i = 12 + 3\sin(\pi t)$ kV, $\phi_i = 0$ and $\zeta_{\rm b} = \zeta_{\rm c} = \frac{2 \pi}{3}$~rad.
\end{itemize}
Figure~\ref{fig:balanced system} shows the phase voltages for E1, E2.  Since the voltage is balanced and the curve in the plane $(\alpha, \beta)$ is a circle, there is a perfect match between the estimations obtained with the geometrical methods, which both return, as expected, a constant frequency ($\omega_a=\omega_{\kappa}= 1$~pu, $\forall t$).  Yet the two methods return the right result for different reasons: $\omega_{\kappa}$ is constant because the circle has a constant curvature; whereas $\omega_a$ is constant because the circle is a special case of an ellipse.  The PLL also works well in E1-E2.

\subsubsection{Unbalanced Voltage}

We consider three examples of unbalanced voltages with constant frequency $\omega_o= 100 \pi $~rad/s: (i) with unequal constant magnitudes; (ii) with unequal and time-varying voltage magnitudes; (iii) with unequal phase displacements.  The following parameters are used:
\begin{itemize}
\item E3: $V_{\rm a} = V_{\rm c} = 12$ kV, $V_{\rm b}= 8$~kV, $\zeta_{\rm b} = \zeta_{\rm c} = \frac{2 \pi}{3}$~rad.
\item E4: $V_{\rm a} = V_{\rm c} = 12 + 3\sin(\pi t)$ kV, $V_{\rm b} = 8 + 2\sin(2 \pi t) $~kV, and $\zeta_{\rm b} = \zeta_{\rm c} = \frac{2 \pi}{3}$ rad.
\item E5: $V_{\rm a} = V_{\rm b}= V_{\rm c} = 12 $~kV, $\zeta_{\rm b} = \frac{-2 \pi}{3}, \zeta_{\rm c} = \frac{1.5 \pi}{3}$~rad.
\end{itemize}

\begin{figure}[ht!]
  \centering
  \begin{subfigure}{0.49\columnwidth}
    \centering
    \includegraphics[width=\linewidth]{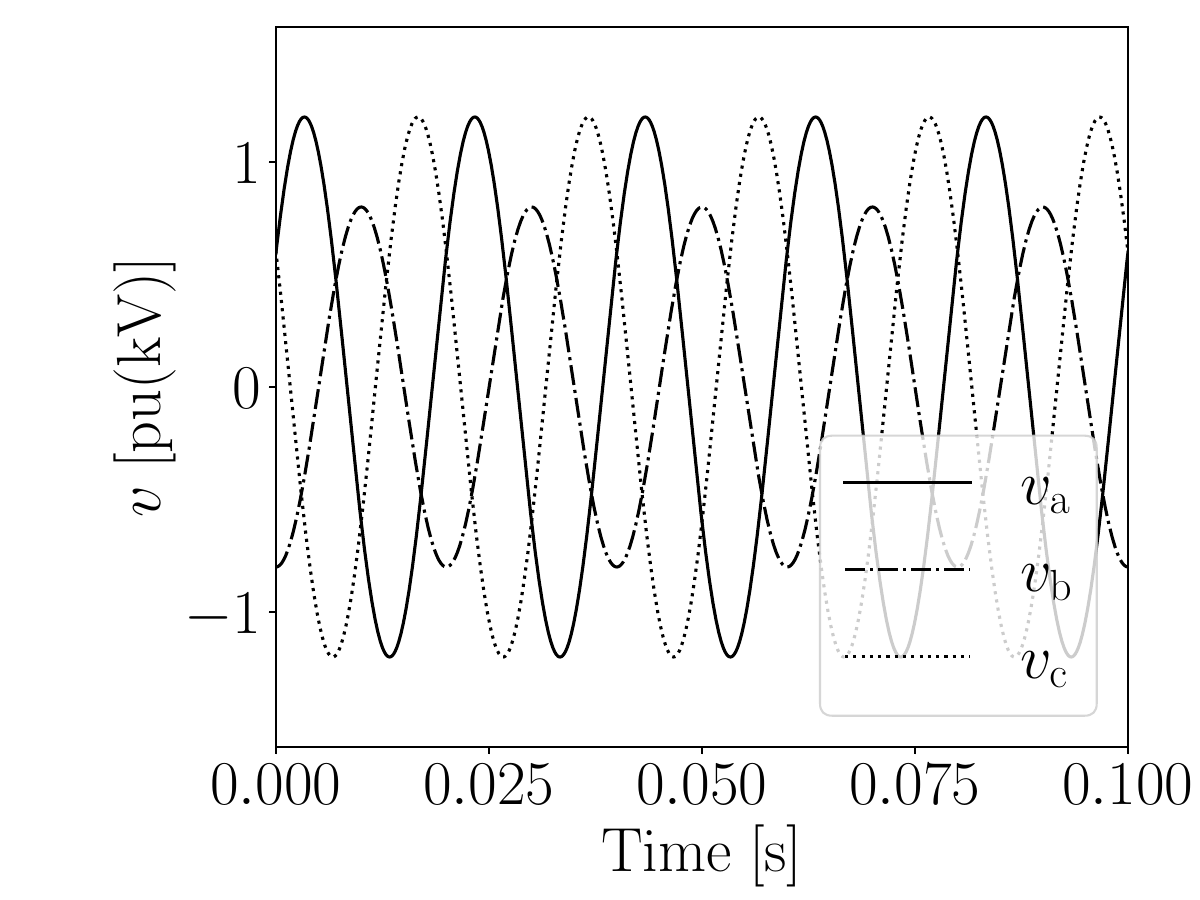}
    \vspace{-6mm}
    \caption{E3: Voltage components}
    \label{fig:E4a}
  \end{subfigure}
  \hfill
  \begin{subfigure}{0.49\columnwidth}
    \centering
    \includegraphics[width=\linewidth]{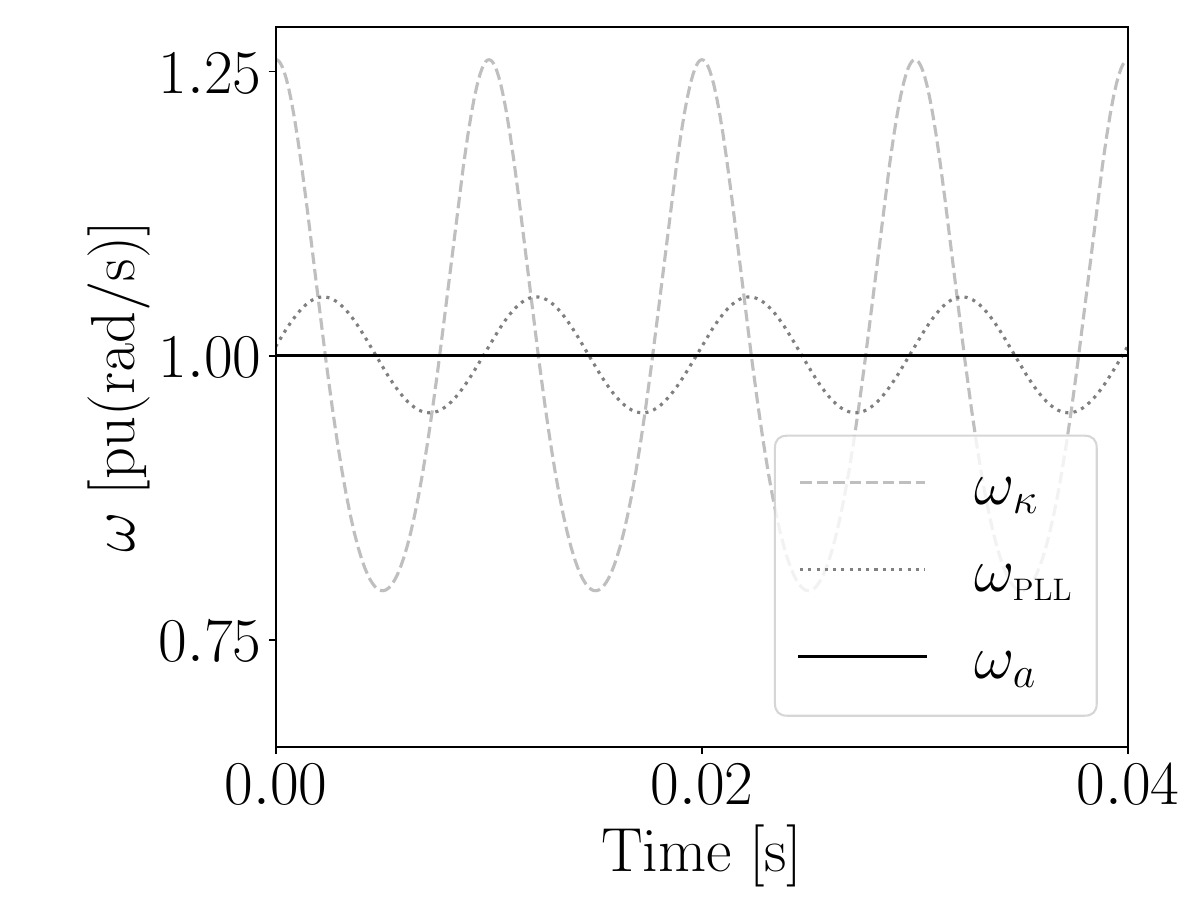}
    \vspace{-6mm}
    \caption{E3: Geometric frequency}
    \label{fig:E4b}
  \end{subfigure}
  \vspace{1mm}
  \begin{subfigure}{0.49\columnwidth}
    \centering
    \includegraphics[width=\linewidth]{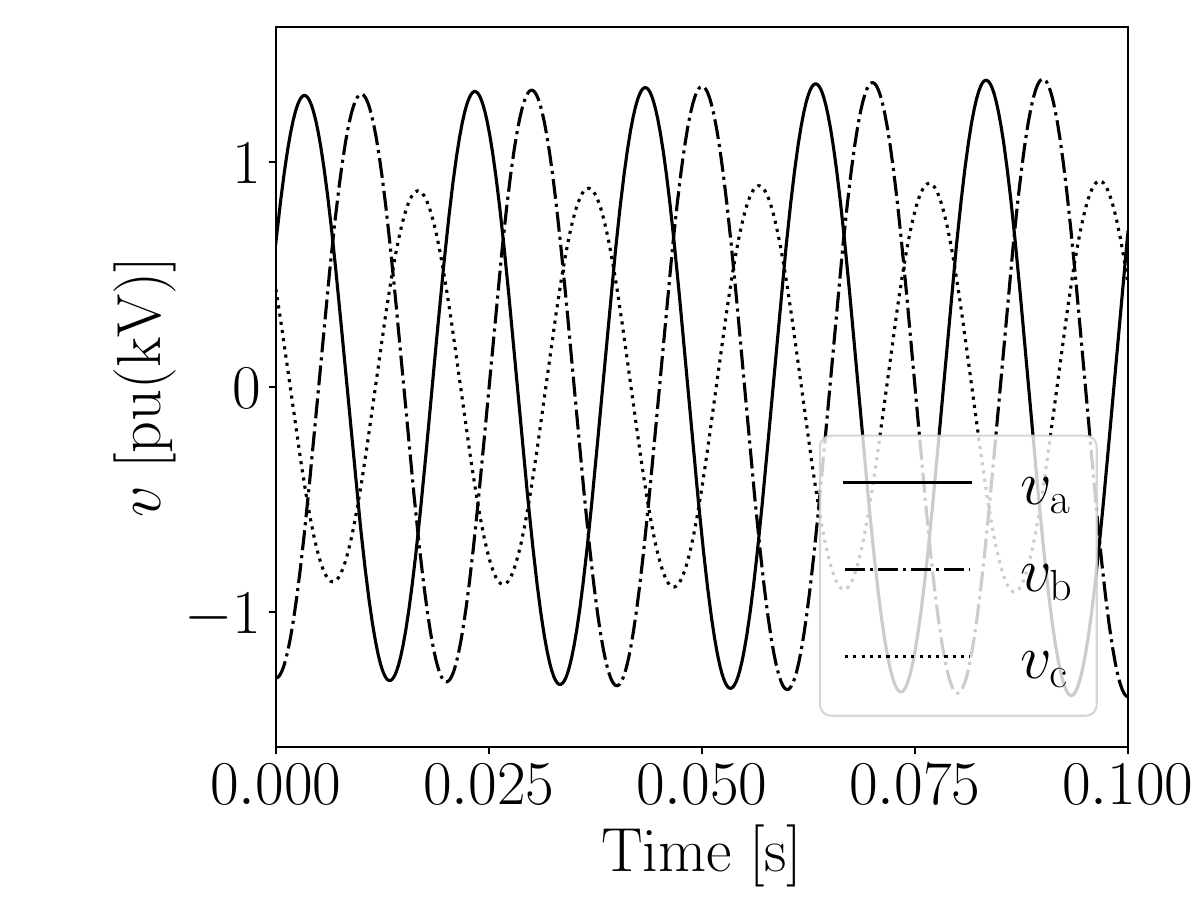}
    \vspace{-6mm}
    \caption{E4: Voltage components}
    \label{case:E5a}
  \end{subfigure}
  \begin{subfigure}{0.49\columnwidth}
    \centering
    \includegraphics[width=\linewidth]{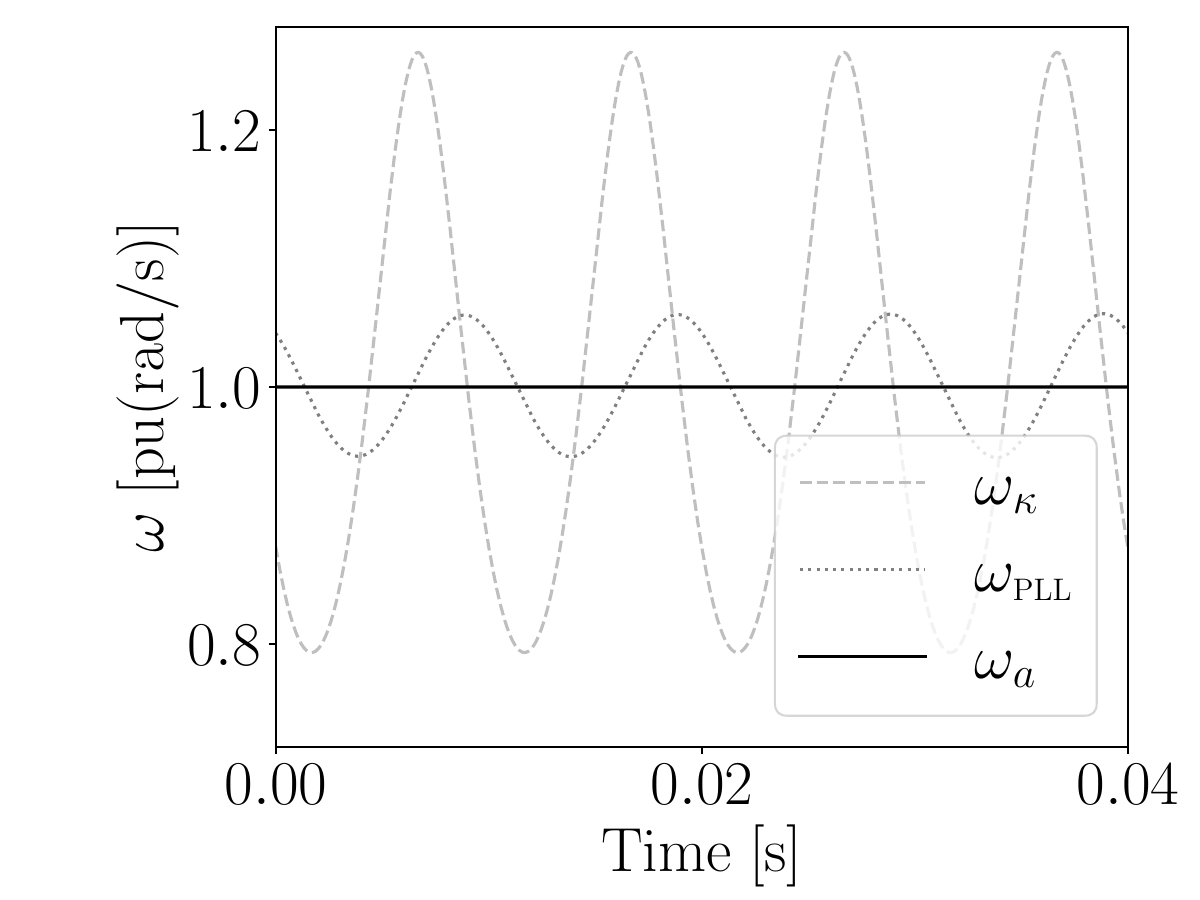}
    \vspace{-6mm}
    \caption{E4: Geometric frequency}
    \label{fig:E5b}
  \end{subfigure}
  \vspace{1mm}
  \begin{subfigure}{0.49\columnwidth}
    \centering
    \includegraphics[width=\linewidth]{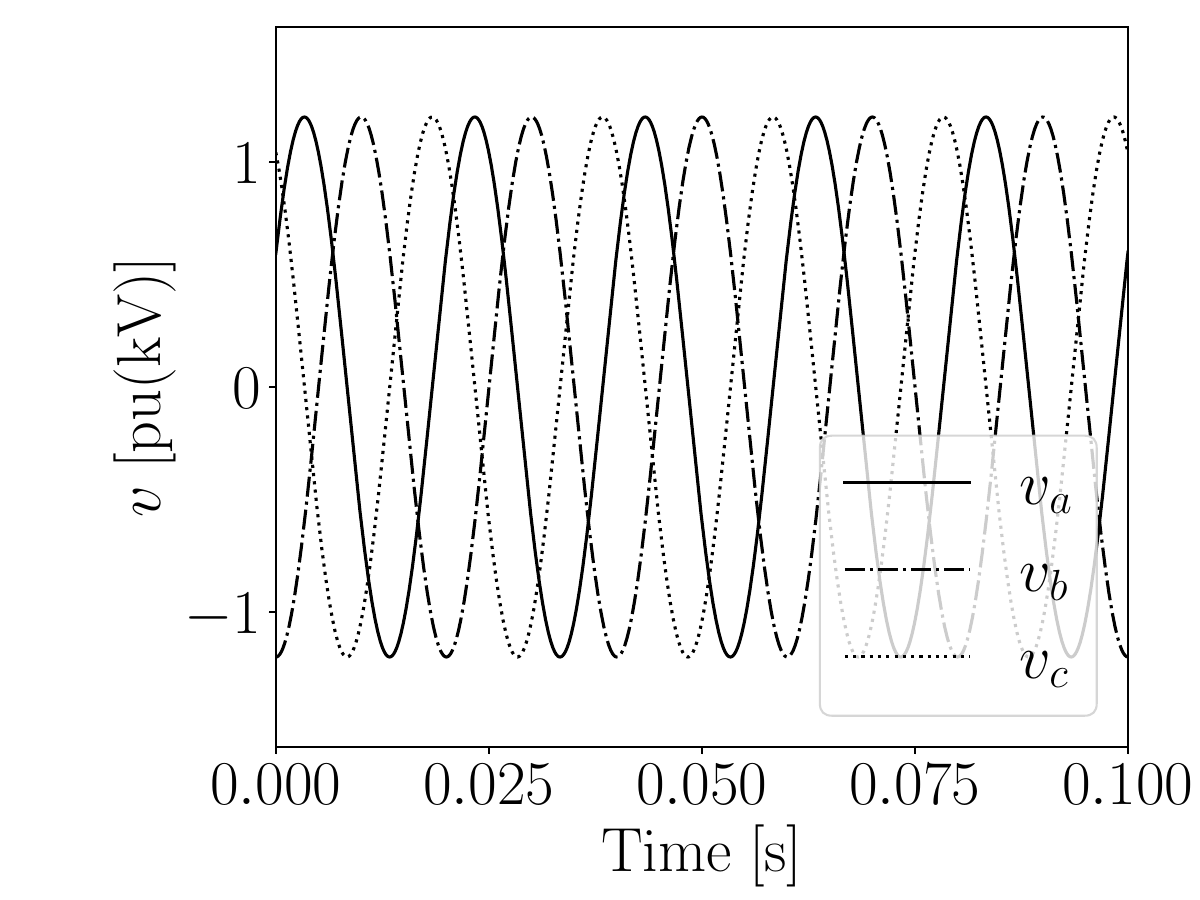}
    \vspace{-6mm}
    \caption{E5: Voltage components}
    \label{fig:E6a}
  \end{subfigure}
  \hfill
  \begin{subfigure}{0.49\columnwidth}
    \centering
    \includegraphics[width=\linewidth]{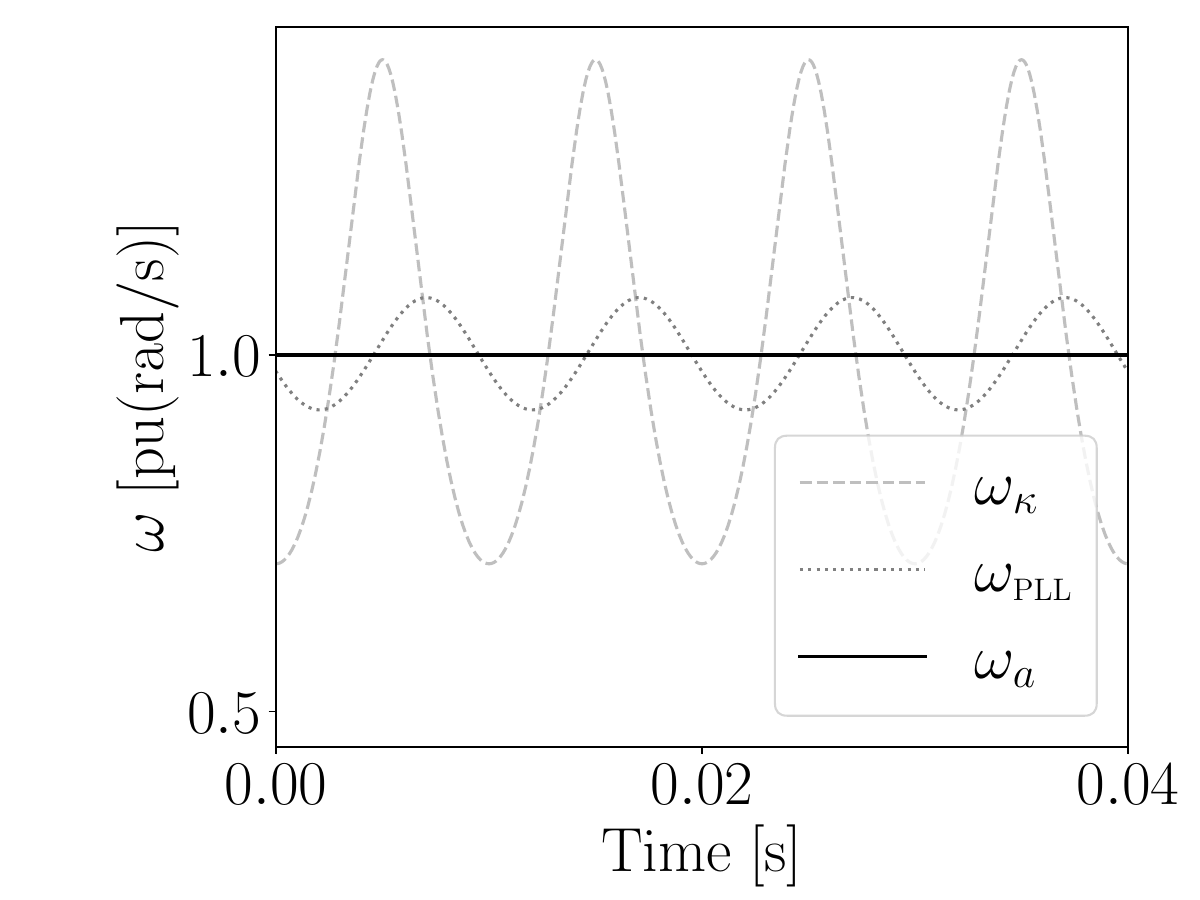}
    \vspace{-6mm}
    \caption{E5: Geometric frequency}
    \label{fig:E6b}
  \end{subfigure}
  \vspace{-5mm}
  \caption{Unbalanced 3-phase voltages and estimated frequencies.}
  \label{fig:unbalanced system}
  \vspace{-2mm}
\end{figure}

For E3-E5, $\phi_i = 0$.  Figure~\ref{fig:unbalanced system} shows the voltage components and estimated geometric and PLL frequencies for E3-E5.  In all three examples, the curves in the $(\alpha, \beta)$ plane are ellipses. This means that the curvature obtained using the Frenet frame is time-varying and periodic, thus leading to a time-varying and periodic $\omega_{\kappa}$.  Moreover, the PLL also outputs a time-varying frequency in the form of a significant ripple around $\omega_o$. On the other hand, \eqref{eq:IF} returns a constant $\omega_a$ equal to $\omega_o$ (in pu), which is consistent with the expected in this case result.

\subsection{Three-Phase Voltage with Time-Varying Frequency}
 
We consider two examples of voltage with varying frequency: (i) a voltage with frequency that varies periodically around its average value.  This case resembles the voltage transient following a contingency in a power system, where voltage phase angle oscillations arising due to electro-mechanical swings of synchronous machines are poorly damped and sustain for several seconds; (ii) an extreme case in power systems, where the voltage components are time-varying and have unequal frequencies.  The following parameters are used:

\begin{itemize}
\item E6: $V_i= 12$~kV, $\zeta_{\rm b} = \zeta_{\rm c}= \frac{2 \pi}{3}$ rad and $\phi_i(t)= \pi \sin (0.4 \pi t)$ rad.
\item E7: $V_i= 12$~kV, $\zeta_{\rm b} = \zeta_{rm c} = \frac{2 \pi}{3}$ rad and $\phi_{\rm a}(t)=\phi_{\rm b}(t)= \pi \sin (0.4 \pi t), \phi_{\rm c}(t) = 1.1 \pi \sin (0.4 \pi t)$ rad.
\end{itemize}

\begin{figure}[ht!]
  \centering
  \begin{subfigure}{0.49\columnwidth}
    \centering
    \includegraphics[width=\linewidth]{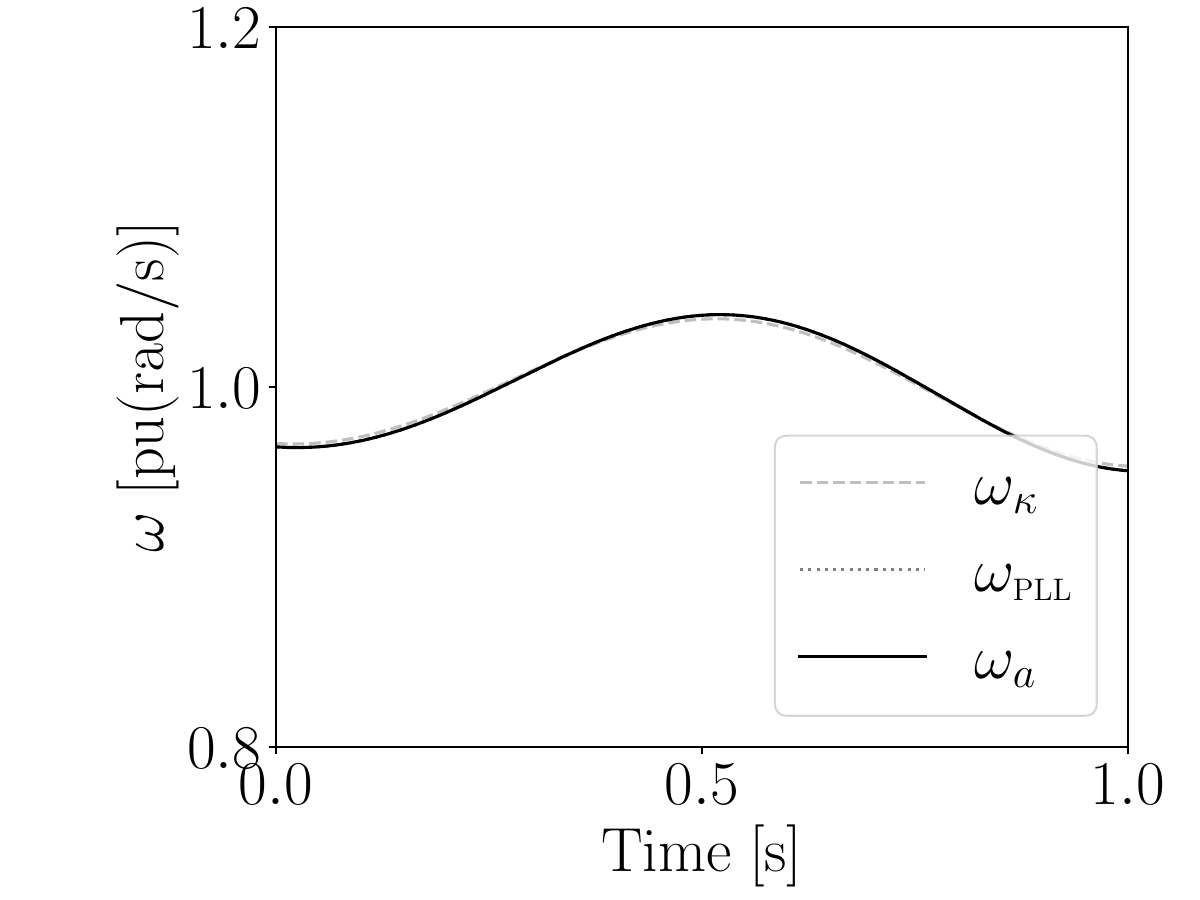}
    \vspace{-6mm}
    \caption{E6}
    \label{fig:E6}
  \end{subfigure}
  \hfill
  \begin{subfigure}{0.49\columnwidth}
    \centering
    \includegraphics[width=\linewidth]{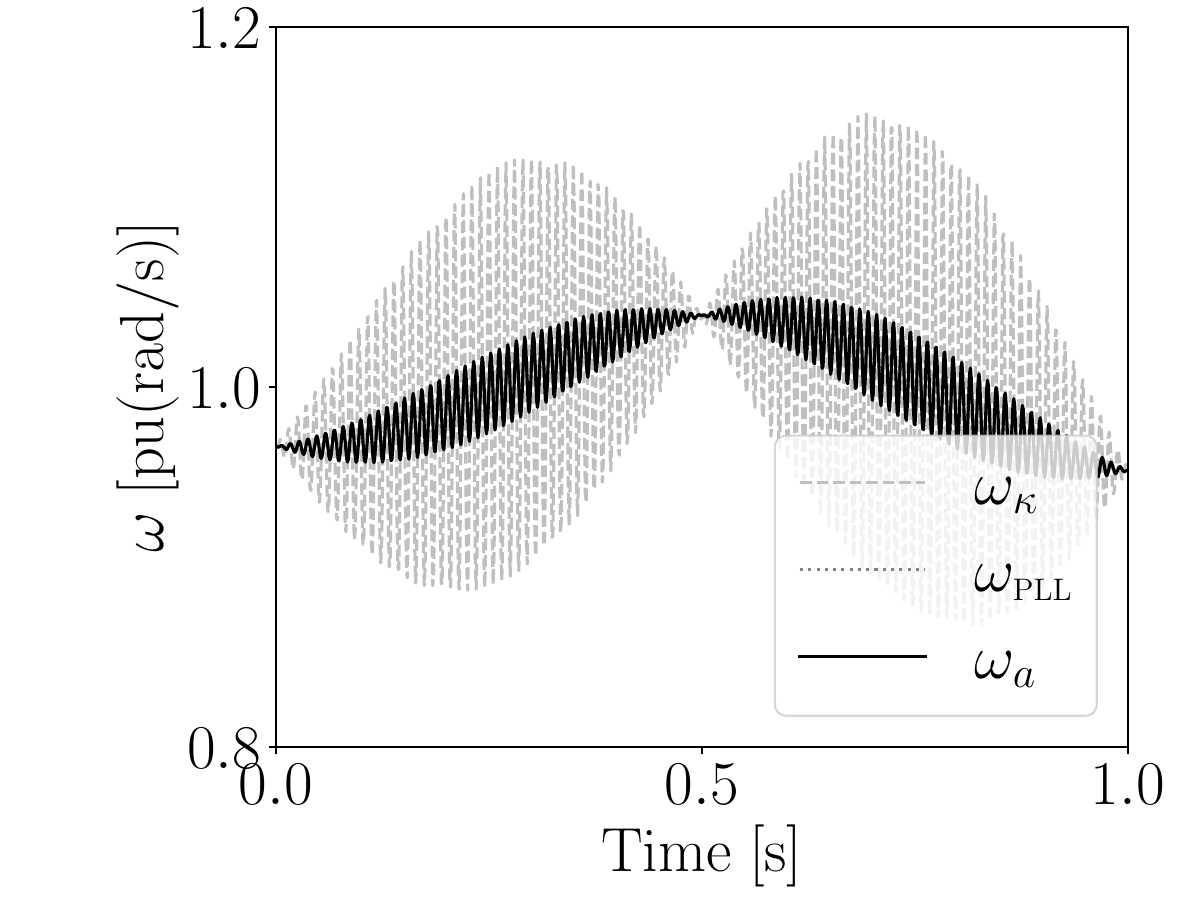}
    \vspace{-6mm}
    \caption{E7}
    \label{fig:E7}
  \end{subfigure}
  \vspace{-6mm}
  \caption{Estimated angular frequency.}
  \label{fig:timevarying}
  \vspace{-4mm}
\end{figure}

Figure~\ref{fig:E6} shows $\omega_a$, $\omega_{\scriptscriptstyle{\rm PLL}}$, and $\omega_\kappa$, for E6.  Despite the approximations imposed by assuming \eqref{eq:cond1}, \eqref{eq:cond2}, we note that \eqref{eq:IF} is able to precisely track the exact instantaneous frequency.  In this example, also $\omega_{\kappa}$ and $\omega_{\scriptscriptstyle{\rm PLL}}$ track well $\rm IF$.  On the other hand, for E7, while $\omega_a$ and $\omega_{\scriptscriptstyle{\rm PLL}}$ still track well the exact frequency, $\omega_{\kappa}$ shows significant fluctuations, see Fig.~\ref{fig:E7}.

\vspace{-5mm}

\vspace{-3mm}

\subsection{Real Voltage Dip Measurements from DFIG in Spain}
\label{sub:realdata}

\begin{figure}[ht!]
  \centering
    \begin{subfigure}{0.49\columnwidth}
    \centering
    \includegraphics[width=\linewidth]{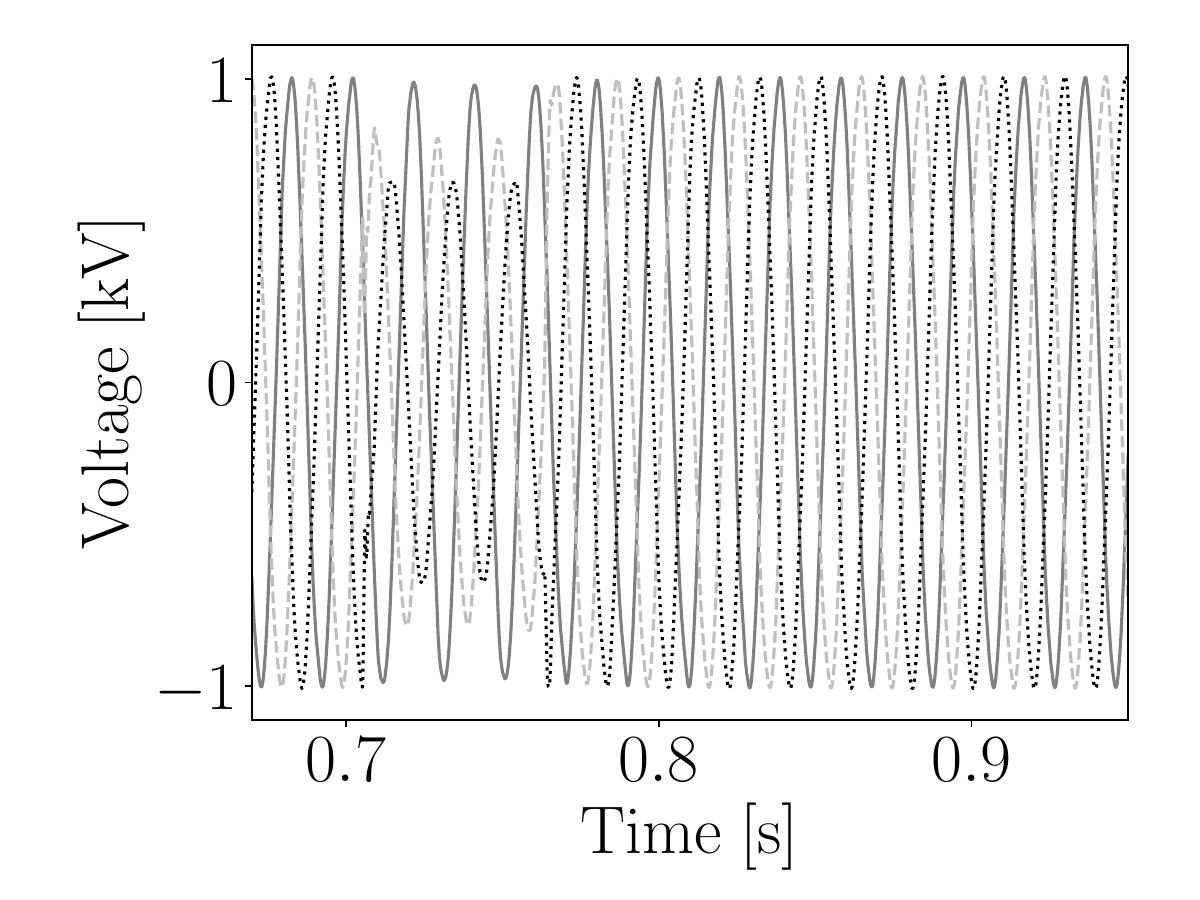}
    \vspace{-7mm}
    \caption{Fault~1: voltage components.}
    \label{fig:c4:v}
  \end{subfigure}
  \hfill
  \begin{subfigure}{0.49\columnwidth}
    \centering
    \includegraphics[width=\linewidth]{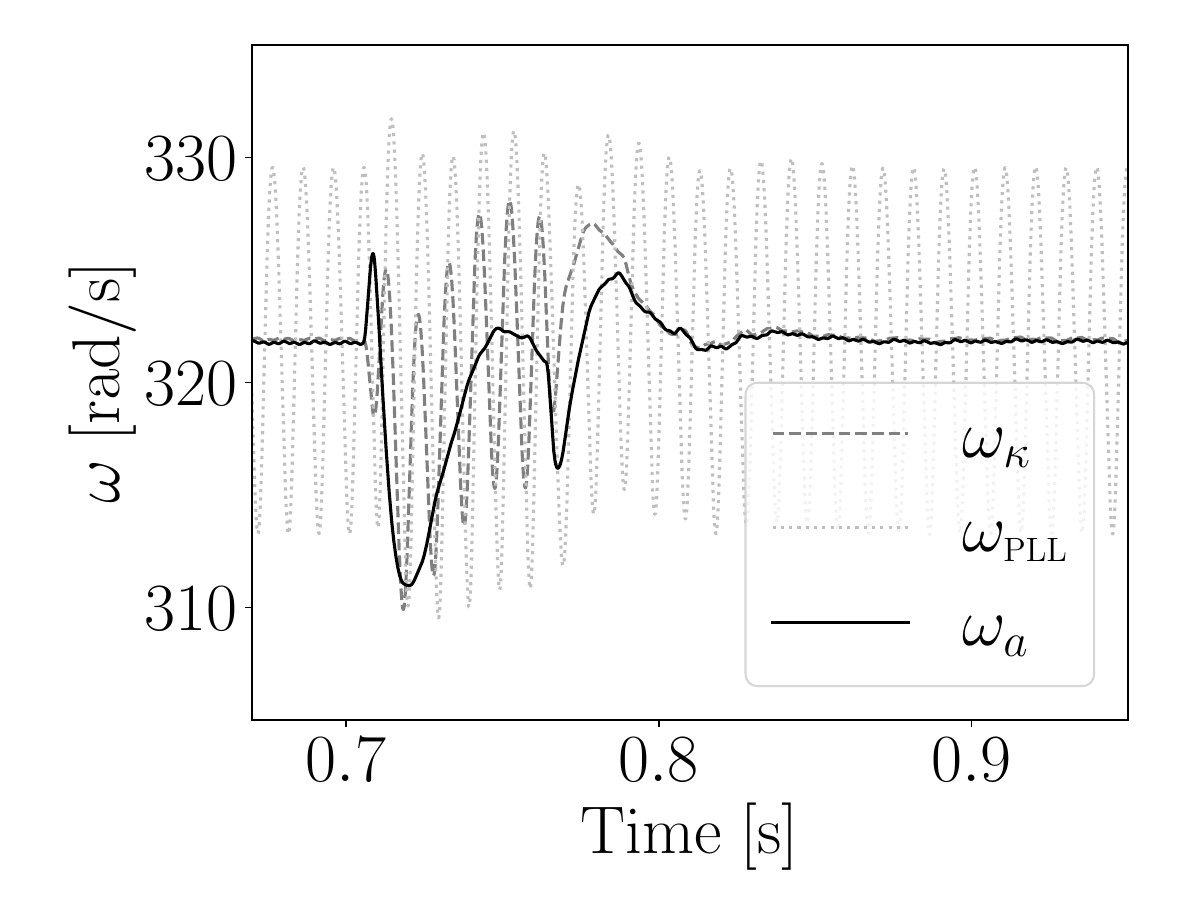}
    \vspace{-7mm}
    \caption{Fault~1: geometric frequency.}
    \label{fig:c4:w}
  \end{subfigure}
    \begin{subfigure}{0.49\columnwidth}
    \centering
    \includegraphics[width=\linewidth]{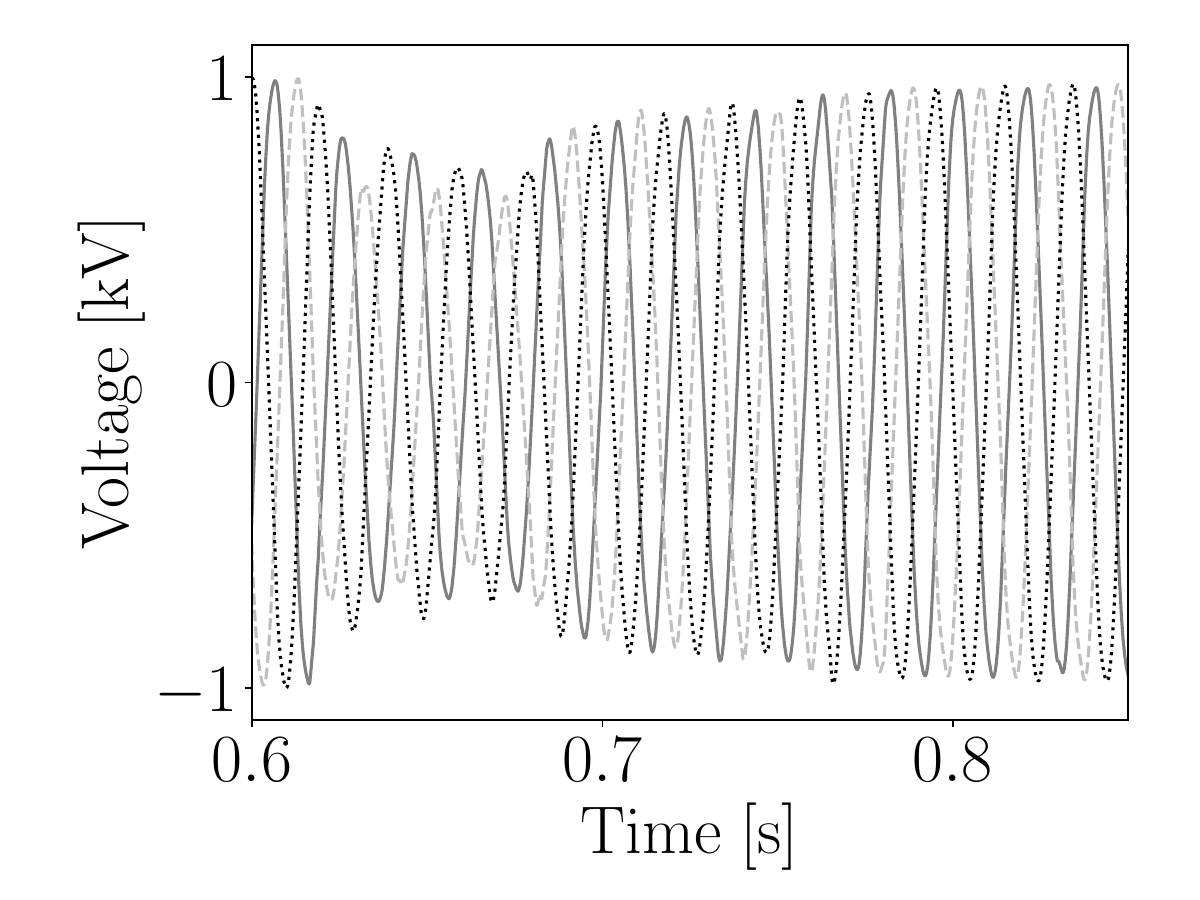}
    \vspace{-7mm}
    \caption{Fault ~2: voltage components.}
    \label{fig:c5:v}
  \end{subfigure}
  \hfill
  \begin{subfigure}{0.49\columnwidth}
    \centering
    \includegraphics[width=\linewidth]{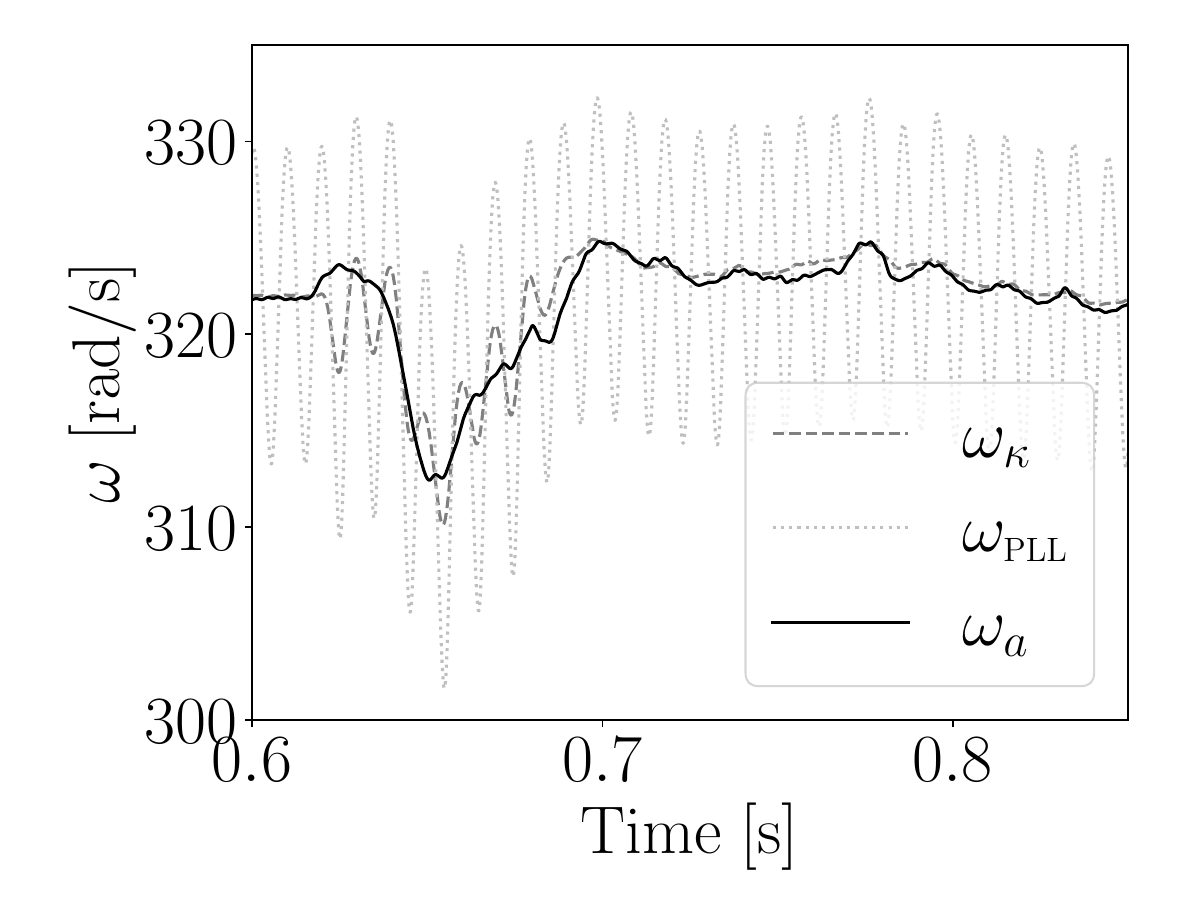}
    \vspace{-7mm}
    \caption{Fault~2: geometric frequency.}
    \label{fig:c5:w}
  \end{subfigure}
  \vspace{-4.5mm}
  \caption{Estimated frequency, real voltage dip data from 690~V DFIG.}
  \label{fig:realdata}
  \vspace{-3mm}
\end{figure}

In this section, the proposed frequency estimation method is further tested using three-phase waveform measurements from two real events of unbalanced faults.  The measurements were taken with a sample rate of 10.25~kHz from the stator of a 690~V, 2.0~MW doubly-fed induction generator~(DFIG) installed in Moralejo, Spain.  Figures~\ref{fig:c4:v} and \ref{fig:c5:v} show the behavior of the three-phase voltages during the two unbalanced voltage dips and following fault clearance.
Figures~\ref{fig:c4:w} and \ref{fig:c5:w} show the results of frequency estimation, indicating that $\omega_a$ is more accurate than $\omega_{\scriptscriptstyle{\rm PLL}}$ and $\omega_\kappa$ for both unbalanced voltage dips.  Note that voltage measurements and frequency estimation outputs for all considered methods are filtered with same second order Butterworth digital filter and an IIR filter.  In these scenarios, $\omega_{\scriptscriptstyle{\rm PLL}}$ shows a bigger ripple than $\omega_{\kappa}$.

\subsection{Single-Phase Voltage}

This example illustrates the performance of \eqref{eq:IF} when applied to a single-phase voltage.  We first consider a voltage with time-varying frequency $\omega t + \phi(t)$ and constant amplitude~$V$:
\begin{equation}
\label{eq:single:ex}
  v(t) = V \sin(\omega_o t + \phi(t)) \, .
\end{equation}
The parameters considered for this example are: $V=12$~kV, $\omega_o=100 \pi$ rad/s and $\phi(t) = 0.05 \omega_o e^{-t} (1-\cos(\pi t))$ rad.

\begin{figure}[ht!]
  \centering
  \begin{subfigure}{0.49\columnwidth}
    \centering
    \includegraphics[width=\linewidth]{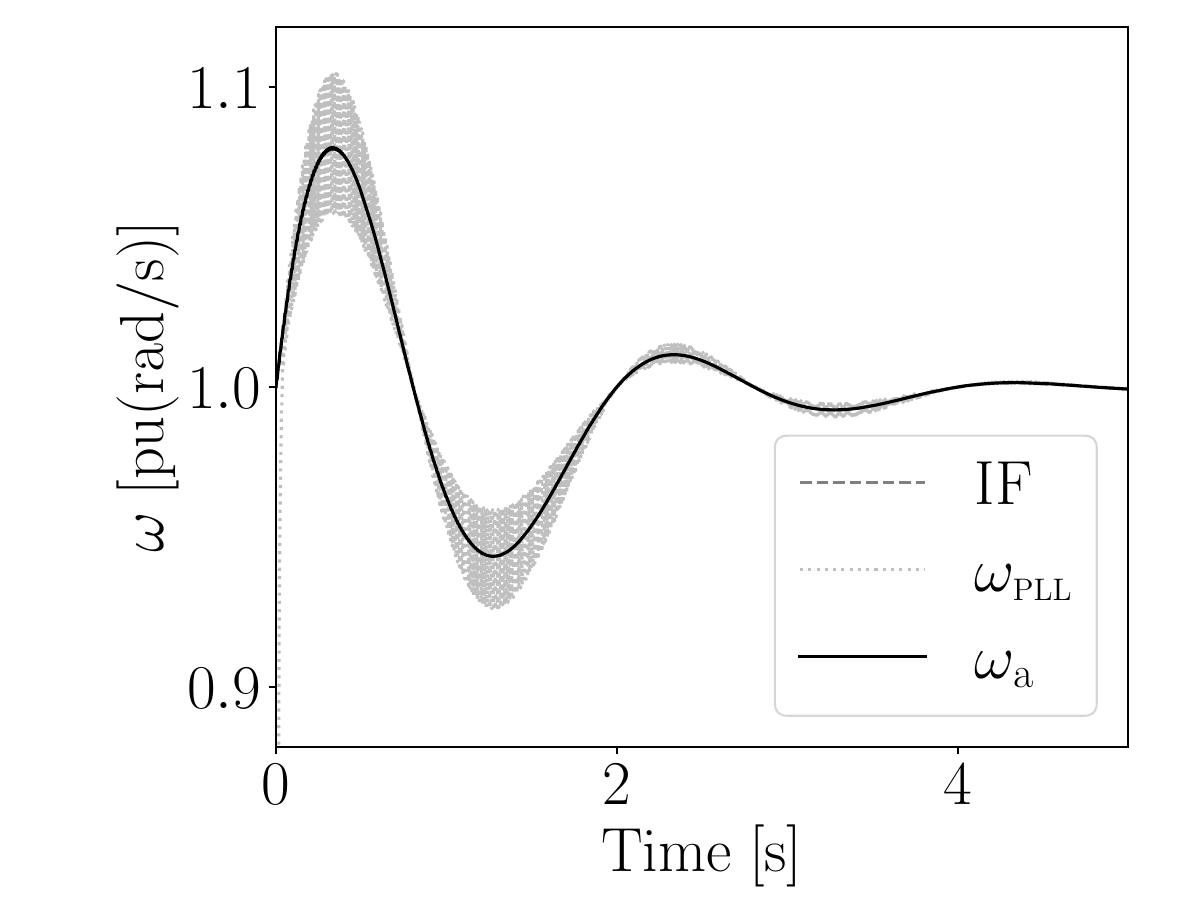}
    \vspace{-6mm}
    \caption{Analytical instantaneous frequency ($\rm IF$) and estimation using a phase shift for the single-phase voltage in \eqref{eq:single:ex}.}
    \label{fig:single phase}
  \end{subfigure}
  \hfill
  \begin{subfigure}{0.49\columnwidth}
    \centering
    \includegraphics[width=\linewidth]{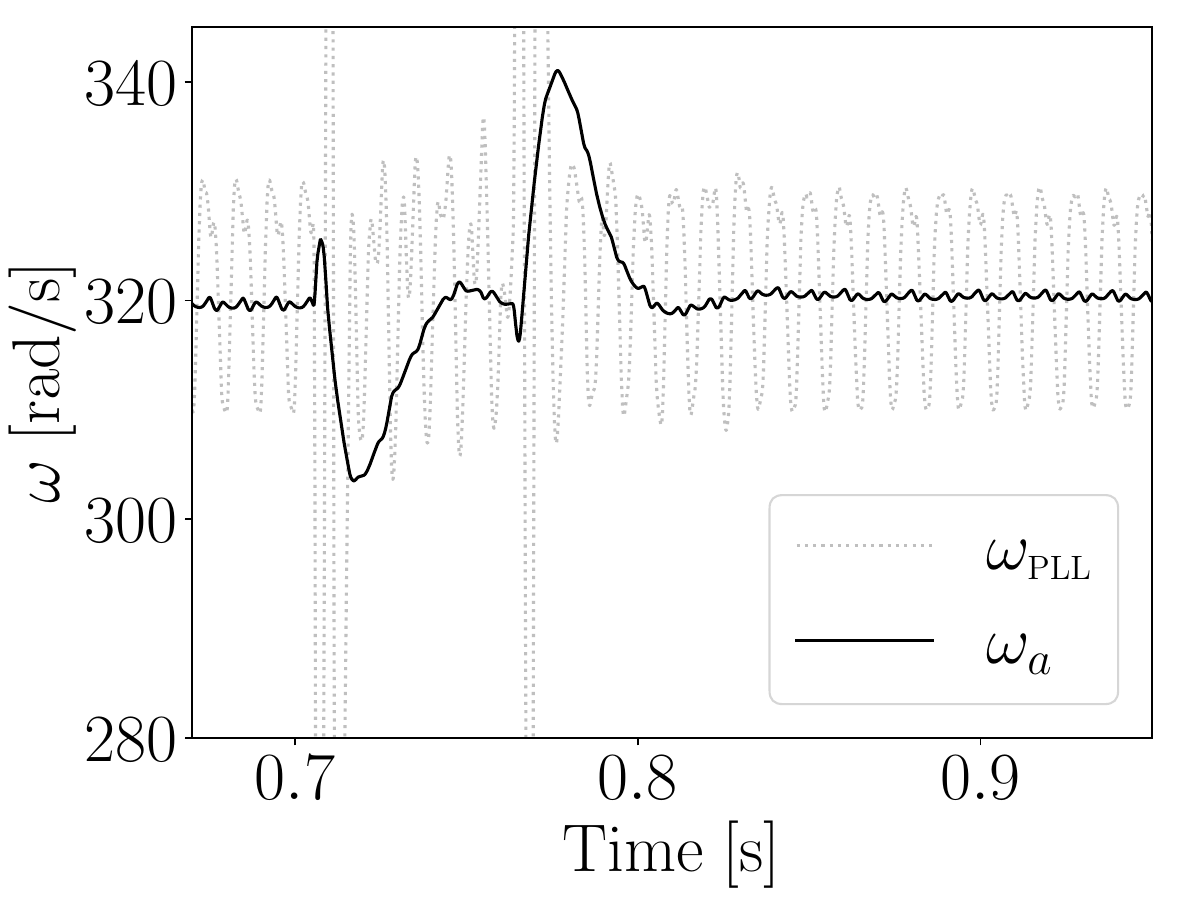}
    \vspace{-6mm}
    \caption{Estimated frequency of $v_{\rm b}$ for real voltage dip data (Fault~1). Comparison of proposed technique with PLL.}
    \label{fig:single phase:realdata}
  \end{subfigure}
  \vspace{-5mm}
  \caption{Estimated angular frequency.}
  \label{fig:single phase_tot}
  \vspace{-6mm}
\end{figure}

As discussed in Section~ref{sub:single}, we construct the second dimension by using the original signal's derivative, as in \eqref{eq:v1ph1}.  Figure~\ref{fig:single phase} illustrates the accuracy of \eqref{eq:IF} in matching the analytical value of the instantaneous frequency $\rm IF = \omega_o + \dot{\phi}$.
Figure~\ref{fig:single phase} also shows the frequency estimated using a conventional PLL where the quadrature signal is obtained using $v(t-\tau)$, where the transport delay is $\tau=0.25 \, T = 0.5\,\pi/\omega_o$.  Despite the approximations resulting from \eqref{eq:cond112}, also in this case $\omega_a$ shows very good accuracy, whereas the PLL shows some ripples due to the fact that the quadrature signal is not exact because the frequency is time-varying.

Finally, we examine the effectiveness of \eqref{eq:IF} in estimating the frequency of a single-phase voltage obtained from real data.  To this end, we use the voltage of phase $\rm b$ ($v_{\rm b}$) from the same voltage dip data considered in Figs.~\ref{fig:c4:v}-\ref{fig:c4:w} (Fault~1).  The results are shown in Fig.~\ref{fig:single phase:realdata}, suggesting that $\omega_a$ provides a better estimation than the PLL.

\section{Conclusions}
\label{sec.Conclusion}

This paper presents a frequency estimation formula based on affine differential geometry for single-phase and unbalanced three-phase voltages.  Approximations based on the nature of typical power system transients are made to achieve a compact explicit expression of the proposed formula.  The adequateness of such approximations is demonstrated through a variety of examples.  When compared to PLLs, as well as to the Frenet frame-based method from \cite{milano2021applications}, the proposed formula proves to be accurate and robust in balanced/unbalanced conditions, as well as for voltages of time-varying magnitude and frequency.


\vfill


\begin{thebibliography}{10}
\providecommand{\url}[1]{#1}
\csname url@samestyle\endcsname
\providecommand{\newblock}{\relax}
\providecommand{\bibinfo}[2]{#2}
\providecommand{\BIBentrySTDinterwordspacing}{\spaceskip=0pt\relax}
\providecommand{\BIBentryALTinterwordstretchfactor}{4}
\providecommand{\BIBentryALTinterwordspacing}{\spaceskip=\fontdimen2\font plus
\BIBentryALTinterwordstretchfactor\fontdimen3\font minus
  \fontdimen4\font\relax}
\providecommand{\BIBforeignlanguage}[2]{{%
\expandafter\ifx\csname l@#1\endcsname\relax
\typeout{** WARNING: IEEEtran.bst: No hyphenation pattern has been}%
\typeout{** loaded for the language `#1'. Using the pattern for}%
\typeout{** the default language instead.}%
\else
\language=\csname l@#1\endcsname
\fi
#2}}
\providecommand{\BIBdecl}{\relax}
\BIBdecl

\bibitem{liu2014three}
B.~Liu, F.~Zhuo, Y.~Zhu, H.~Yi, and F.~Wang, ``A three-phase {PLL} algorithm
  based on signal reforming under distorted grid conditions,'' \emph{IEEE
  Trans.~on Power Electronics}, vol.~30, no.~9, pp. 5272--5283, 2014.

\bibitem{song2022fast}
J.~Song, A.~Mingotti, J.~Zhang, L.~Peretto, and H.~Wen, ``Fast
  iterative-interpolated {DFT} phasor estimator considering out-of-band
  interference,'' \emph{IEEE Trans. on Instr. and Meas.}, vol.~71, pp. 1--14,
  2022.

\bibitem{reza2016accurate}
S.~Reza, M.~Ciobotaru, and V.~G. Agelidis, ``Accurate estimation of
  single-phase grid voltage fundamental amplitude and frequency by using a
  frequency adaptive linear {K}alman filter,'' \emph{IEEE J. of Emerging and
  Selected Topics in Power Elec.}, vol.~4, no.~4, pp. 1226--1235, 2016.

\bibitem{nie2019detection}
X.~Nie, ``Detection of grid voltage fundamental and harmonic components using
  {K}alman filter based on dynamic tracking model,'' \emph{IEEE Trans. on Ind.
  Electronics}, vol.~67, no.~2, pp. 1191--1200, 2019.

\bibitem{pradhan2005freq}
A.~Pradhan, A.~Routray, and A.~Basak, ``Power system frequency estimation using
  least mean square technique,'' \emph{IEEE Trans. on Power Delivery}, vol.~20,
  no.~3, pp. 1812--1816, 2005.

\bibitem{hadjidemetriou2016synchronization}
L.~Hadjidemetriou, Y.~Yang, E.~Kyriakides, and F.~Blaabjerg, ``A
  synchronization scheme for single-phase grid-tied inverters under harmonic
  distortion and grid disturbances,'' \emph{IEEE Trans.~on Power Electronics},
  vol.~32, no.~4, pp. 2784--2793, 2016.

\bibitem{santos2008comparison}
S.~Filho \emph{et~al.}, ``Comparison of three single-phase {PLL} algorithms for
  ups applications,'' \emph{IEEE Trans. on Ind. Electronics}, vol.~55, no.~8,
  pp. 2923--2932, 2008.

\bibitem{hao2007measuring}
P.~Hao, W.~Zanji, and C.~Jianye, ``A measuring method of the single-phase ac
  frequency, phase, and reactive power based on the hilbert filtering,''
  \emph{IEEE Trans. on Instr. and Meas.}, vol.~56, no.~3, pp. 918--923, 2007.

\bibitem{sahoo2021phase}
A.~Sahoo, J.~Ravishankar, and C.~Jones, ``Phase-locked loop independent
  second-order generalized integrator for single-phase grid synchronization,''
  \emph{IEEE Trans. on Instr. and Meas.}, vol.~70, pp. 1--9, 2021.

\bibitem{7562505}
F.~Xiao, L.~Dong, L.~Li, and X.~Liao, ``A frequency-fixed {SOGI}-based {PLL}
  for single-phase grid-connected converters,'' \emph{IEEE Trans. on Power
  Electronics}, vol.~32, no.~3, pp. 1713--1719, 2017.

\bibitem{reza2019three}
M.~S. Reza \emph{et~al.}, ``Three-phase {PLL} for grid-connected power
  converters under both amplitude and phase unbalanced conditions,'' \emph{IEEE
  Trans. on Ind. Electronics}, vol.~66, no.~11, pp. 8881--8891, 2019.

\bibitem{karimi2004estimation}
H.~Karimi, M.~Karimi-Ghartemani, and M.~R. Iravani, ``Estimation of frequency
  and its rate of change for applications in power systems,'' \emph{IEEE Trans.
  on Power Delivery}, vol.~19, no.~2, pp. 472--480, 2004.

\bibitem{meral2012improved}
M.~E. Meral, ``Improved phase-locked loop for robust and fast tracking of three
  phases under unbalanced electric grid conditions,'' \emph{IET Generation,
  Transmission \& Distribution}, vol.~6, no.~2, pp. 152--160, 2012.

\bibitem{kulkarni2015design}
A.~Kulkarni and V.~John, ``Design of synchronous reference frame phase-locked
  loop with the presence of dc offsets in the input voltage,'' \emph{IET Power
  Electronics}, vol.~8, no.~12, pp. 2435--2443, 2015.

\bibitem{golestan2012design}
S.~Golestan, M.~Monfared, and F.~D. Freijedo, ``Design-oriented study of
  advanced synchronous reference frame phase-locked loops,'' \emph{IEEE
  Trans.~on Power Electronics}, vol.~28, no.~2, pp. 765--778, 2012.

\bibitem{escobar2014cascade}
G.~Escobar \emph{et~al.}, ``Cascade three-phase {PLL} for unbalance and
  harmonic distortion operation {(CSRF-PLL)},'' in \emph{IECON}, 2014, pp.
  5489--5493.

\bibitem{milano2021applications}
F.~Milano, G.~Tzounas, I.~Dassios, and T.~Kerci, ``Applications of the {Frenet}
  frame to electric circuits,'' \emph{IEEE Trans. on Circuits and Systems I:
  Regular Papers}, vol.~69, no.~4, pp. 1668--1680, 2021.

\bibitem{milano2022paradox}
F.~Milano, G.~Tzounas, I.~Dassios, M.~A.~A. Murad, and T.~K{\"e}r{\c{c}}i,
  ``Using differential geometry to revisit the paradoxes of the instantaneous
  frequency,'' \emph{IEEE Open Access J. of Power and Energy}, 2022.

\bibitem{milano2021geometrical}
F.~Milano, ``A geometrical interpretation of frequency,'' \emph{IEEE Trans. on
  Power Systems}, vol.~37, no.~1, pp. 816--819, 2021.

\bibitem{milano2022frenet}
F.~Milano{}, ``The {Frenet} frame as a generalization of the {Park}
  transform,'' \emph{IEEE Trans. on Circuits and Systems I: Regular Papers},
  vol.~70, no.~2, pp. 966--976, 2022.

\bibitem{lewis2018bountiful}
A.~D. Lewis, ``The bountiful intersection of differential geometry, mechanics,
  and control theory,'' \emph{Annual Review of Control, Robotics, and
  Autonomous Systems}, vol.~1, pp. 135--158, 2018.

\bibitem{craizer2006parabolic}
M.~Craizer, T.~Lewiner, and J.-M. Morvan, ``Parabolic polygons and discrete
  affine geometry,'' in \emph{Brazilian Symposium on Computer Graphics and
  Image Processing}.\hskip 1em plus 0.5em minus 0.4em\relax IEEE, 2006, pp.
  19--26.

\bibitem{flash2007affine}
T.~Flash and A.~Handzel, ``Affine differential geometry analysis of human arm
  movements,'' \emph{Biological Cybernetics}, vol.~96, pp. 577--601, 2007.

\bibitem{nomizu1994affine}
K.~Nomizu and T.~Sasaki, \emph{Affine differential geometry: geometry of affine
  immersions}.\hskip 1em plus 0.5em minus 0.4em\relax Cambridge University
  Press, 1994.

\bibitem{calabi1996affine}
E.~Calabi, P.~J. Olver, and A.~Tannenbaum, ``Affine geometry, curve flows, and
  invariant numerical approximations,'' \emph{Advances in Mathematics}, vol.
  124, no.~1, pp. 154--196, 1996.

\end{thebibliography}
\end{document}